\documentclass{svproc}
\usepackage{url}

\usepackage{amsmath} 
\usepackage{amssymb}  
\usepackage[ruled,longend]{algorithm2e}
\usepackage{multicol}
\usepackage{graphicx}
\usepackage{caption}
\usepackage{threeparttable}
\usepackage{subcaption}
\usepackage{cite}
\usepackage{float}
\usepackage{color}
\begin{document}
\mainmatter              
\title{D2C 2.0:  Decoupled Data-Based Approach for Learning to Control Stochastic Nonlinear Systems via Model-Free ILQR}
\titlerunning{D2C-2.0}  
%
\author{Karthikeya S Parunandi$^1$, Aayushman Sharma$^2$, Suman Chakravorty$^3$ and Dileep Kalathil$^4$}

\authorrunning{Karthikeya Parunandi et al.} 
%
\tocauthor{Karthikeya S Parunandi$^1$, Aayushman Sharma$^2$, Suman Chakravorty$^3$ and Dileep Kalathil$^4$}
\institute{Texas A$\&$M University, College Station TX 77843, USA\\
\email{karthikeyasharma91@gmail.com}$^1$, \email{aayushmansharma@tamu.edu$^2$}, \email{s.chakrav@tamu.edu$^3$}, and \email{dileep.kalathil@tamu.edu$^4$}
}

\maketitle              

\begin{abstract}
In this paper, we propose a structured linear parameterization of a feedback policy to solve the model-free stochastic optimal control problem. This parametrization is corroborated by a decoupling principle that is shown to be near-optimal under a small noise assumption, both in theory and by empirical analyses. Further, we incorporate a model-free version of the Iterative Linear Quadratic Regulator (ILQR) in a sample-efficient manner into our framework. Simulations on systems over a range of complexities reveal that the resulting algorithm is able to harness the superior second-order convergence properties of ILQR. As a result, it is fast and is scalable to a wide variety of higher dimensional systems. Comparisons are made with a state-of-the-art reinforcement learning algorithm, the Deep Deterministic Policy Gradient (DDPG) technique, in order to demonstrate the significant merits of our approach in terms of training-efficiency.
\end{abstract}
\section{INTRODUCTION}
\label{introduction}
It is well-known that a general solution to stochastic optimal control problem is computationally intractable. Hence, most methods rely on approximated dynamic programming methods (See, for e.g., \cite{approximate_dp}). In case of an unknown system model, the problem is further complicated. This problem of learning to optimally control an unknown general nonlinear system has been well-studied in the stochastic adaptive control literature \cite{Ioannou} \cite{Astrom}. A very similar problem is also addressed in the field of reinforcement learning (RL) except that it specifically aims to solve a discrete-time stochastic control process (Markov Decision Process), more often with an infinite horizon formulation \cite{Sutton}. Also, RL has attracted a great attention with its recent success, especially in its sub-domain of deep RL \cite{atari_nature} \cite{Silver} \cite{Lillicrap}.\\
Most methods in the RL can be divided into model-based and model-free approaches, based on whether a model is constructed before solving for an optimal policy or not, respectively. Among them, model-based RL typically requires lesser number of samples and can be global or local model-based. Global models are typically approximated with Gaussian Processes (GPs) or Neural Networks, and hence can be inaccurate in some portions of its state space. Local models (such as in \cite{Kumar2016}\cite{Mitrovic2010}) in contrast offer more sample efficient and accurate solutions. On the other hand, model-free RL, since has to bypass learning the system model requires substantially more number of training samples to train the policy \cite{inaccurate models in rl}. While an overwhelming number of successful works in the recent past have relied on a complex parameterizations such as by using neural networks \cite{Lillicrap} \cite{acktr}\cite{trpo} \cite{ppo}, they are challenged with non-determinism in reproducing the results and in benchmarking against existing methods \cite{Henderson}\cite{Reproducibility_RL}. Recent works such as by \cite{Rajeswaran} showed that a simple linear parameterization can be very effective and can result in performances competitive to that of the former methods. The current work, aimed at deriving a simple but structured, sample-efficient and thereby a fast policy, falls under the local model-based regime, though we won't explicitly build the system model but rather deal with its gradients in the neighborhood of a nominal trajectory.\\

In the prequel of this work, called D2C \cite{d2c}, instead of solving for the optimal closed-loop policy via dynamic programming as in typical RL techniques, our approach searches for a computationally tractable policy that is nonetheless near-optimal to second order in a noise parameter $\epsilon$. This is done via decoupling of the open and the closed loops, wherein we first solve for an optimal open-loop trajectory, and then design an optimal linear feedback law to stabilize the linear perturbations about the nominal trajectory. The success of this approach was previously demonstrated in \cite{Rafieisakhaei}\cite{Karthikeya} for model-based problems. This way of decoupling results in a compact parameterization of the feedback law in terms of an open loop control sequence and a linear feedback policy, and thereby, constitutes a novel approach to the RL problem via the decoupled search.
\textbf{The contributions of this sequel paper are as follows:}\\
i) We incorporate a data-efficient, and hence, significantly faster way to compute the gradients in the ILQR framework in a model free fashion. This, in turn, is incorporated in our D2C framework to solve for the optimal open-loop trajectory. This is in contrast to the `direct' gradient-descent based approach introduced in the earlier version of our algorithm (D2C 1.0).\\
ii) We make a comparison with a state-of-the-art reinforcement learning algorithm - DDPG  \cite{Lillicrap} as well as D2C 1.0, and show the significant merit of the current method, termed D2C 2.0, in terms of its computational efficiency.\\
iii) D2C 2.0  also constitutes a novel approach to solving complex robotic motion planning problems under uncertainty in a data efficient model-free fashion.

\section{RELATED WORK}
\label{related_work}
Differential Dynamic Programming (DDP) \cite{mayne's ddp} is a class of iterative algorithms for trajectory optimization. One salient aspect of DDP is that it exhibits quadratic convergence for a general non-linear optimal control problem and can be solved in a single step for linear systems. In spite of its origins in the 1960s, it gained popularity only in the last decade due to the success of the modified algorithm called ILQR \cite{ilqr-1}. Though DDP (theoretically) guarantees a quadratic convergence, some of the terms in it involve computing second order derivatives of the system dynamical models. Since the dynamical models of most systems are multivariable vector-valued functions and their second order derivatives being third order tensors, DDP in its original form was not effective for practical implementation. ILQR \cite{ilqr} \cite{ILQG} dropped these terms and introduced regularization schemes \cite{ddp with nonlinear constraints}\cite{ILQG_complex_behaviors} to handle the absence of these computationally expensive terms. This resulted in a faster and stable convergence to the optimal solution.

One advantage that ILQR has in terms of solving it in a model-free manner is that this algorithm is explicit in terms of system model and their first order derivatives. Earlier works such as by \cite{ILQG_complex_behaviors} and \cite{Levine} employed finite differencing in computing the Jacobians in ILQR. Typically, a forward Euler scheme is chosen to independently determine each element ({\it i.e,} gradient) in a Jacobian matrix.  \cite{Mitrovic2010} presented ILQR-LD by learning the dynamics through Locally Weighted Projection Regression (LWPR) approximation of the model first and then obtaining a corresponding analytical Jacobian. Their work also demonstrated a better computational efficiency over finite differences.

Spall introduced the `Simultaneous Perturbation Stochastic Approximation' (SPSA) \cite{spsa} method that evaluates the function only twice to calculate the Jacobian of the cost function. In this paper, a similar formulation is derived to compute the Jacobians of the system model online, in a model-free fashion, through the least squares method in a central-difference scheme.

The rest of the paper is organized as follows. Section \ref{problem_formulation} introduces the terminology and provides the problem statement. In section \ref{decoupling_principle}, we present a near-optimal decoupling principle that forms the basis of our algorithm. Section \ref{algorithm} describes our algorithm - D2C 2.0 in detail. In Sections \ref{simulations} and \ref{discussion}, we provide results and a discussion of the application of the D2C 2.0 to several benchmark examples, along with comparisons to DDPG and the D2C technique.
\section{PROBLEM FORMULATION}\label{problem_formulation}
Let ${\bf x_t}$ $\in$ $\mathcal{X} \in$  $\mathbb{R}^{n_x}$ represent the state of a system and ${\bf u_t} \in \mathcal{X} \in \mathbb{R}^{n_u}$ be the control signal at time $t$ respectively. 
Let $f : \mathbb{R}^{n_x} \times \mathbb{R}^{n_u} \rightarrow \mathbb{R}^{n_x}$ denote the deterministic state transition model of the system. Let $\pi_t : \mathbb{R}^{n_x} \rightarrow \mathbb{R}^{n_u}$ be the policy to be applied on the system. Let $c_t : \mathbb{R}^{n_x} \times \mathbb{R}^{n_u} \rightarrow \mathbb{R}$ be the time-indexed incremental cost function for all $t\in\{0,1,\cdots,N-1\}$ and $c_N : \mathbb{R}^{n_x} \rightarrow \mathbb{R}$ be the terminal cost function. Let $J_t^{\pi_t} : \mathbb{R}^{n_x} \rightarrow \mathbb{R}$ represent the cost-to-go function under the policy $\pi_t$ and $Q_t^{\pi_t}: \mathbb{R}^{n_x} \times \mathbb{R}^{n_u} \rightarrow \mathbb{R}$ be the corresponding action-value function, both at time $t$. 
Finally, let us denote the derivatives of a variable at time $t$ w.r.t. the state ${\bf x_t}$ or the control ${\bf u_t}$ by placing respective variables in the subscript. For example, $J_{\bf x_t} = \frac{\partial J_t}{\partial {\bf x_t}}|_{\bf \bar{x}_t}$. Unless stated otherwise, they are evaluated at the nominal state or control at time $t$. 

Now, let us define the notion of process noise in the system through an important parameter called `noise scaling factor', denoted by $\epsilon$. We assume that the noise is additive and white, zero-mean Gaussian. Hence, the state evolution equation of our stochastic system is represented as ${\bf x_{t+1}} = f({\bf x_t}, {\bf u_t}) + \epsilon \omega_t$, where $\epsilon < 1$ and $\omega_t$ is zero-mean Gaussian distributed.\\
Let us define the notion of a `nominal' trajectory, given a policy $\pi_t(.)$, as follows - the nominal control actions of the policy are the control actions output by the policy when all the noise variables $\omega_t =0$. Let us represent the nominal trajectory variables with \textit{bars} over them $i.e,$ ${\bf \bar{x}_t}$ and ${\bf \bar{u}_t}$, both at time $t$. So, the nominal trajectory is given by, $\mathbb{T}_{nom} = \{\bf \bar{x}_{0:N}, \bar{u}_{0:N-1}\}$, where $N$ is the time horizon, and $\bar{\mathbf{x}}_{t+1} = f(\bar{\mathbf{x}_t}, \bar{\mathbf{u}}_t)$, where $\bar{\mathbf{u}}_t = \pi_t(\bar{\mathbf{x}}_t)$. 


Using the above definitions, the problem of {\it stochastic optimal control} that we consider is to find the set of time-varying optimal control policies $\pi = \{\pi_1, \pi_2, \cdots, \pi_{N-1}\}$ (subscripts denote the time-indices) that minimize the finite horizon total expected cumulative cost, $i.e,$
\begin{equation}
    \tilde{J}^{\pi} = \mathbb{E}_{\pi}[J^{\pi}] = \mathbb{E}_{\pi}[\sum_{t=0}^{N-1}c_t({\bf x_t}, \pi_t({\bf x_t})) + c_N({\bf x_N})],
\label{soc_problem}
\end{equation}
such that the system model {\it i.e,} ${\bf x_{t+1}} = f({\bf x_t}, \pi_t({\bf x_t})) + \epsilon \omega_t$, relevant constraints and boundary conditions are satisfied. Since the current approach deals with model-free problem, the system model and its constraints are assumed to be implicitly satisfied in the data that the system/agent gathers from its environment (through a simulator or a real world experiment). 
\section{A NEAR-OPTIMAL DECOUPLING PRINCIPLE}
\label{decoupling_principle}
{\textit{Note : The results provided in this particular section are restated from our prequel paper on D2C primarily for the reader's convenience and for the sake of completeness.}}\\
In this section, we describe a near-optimal decoupling principle which provides the theoretical foundation for our D2C technique.\\
Assuming that $f(.)$ and $\pi_t(.)$ are sufficiently smooth, let us describe the linearized dynamics about the nominal trajectory $\mathbb{T}_{nom}$ corresponding to the policy $\pi_t(.)$ as follows. The state and control perturbations, given by $\delta {\bf x_t} = {\bf x_t} - {\bf \bar{x}_t}$ and $\delta {\bf u_t} = {\bf u_t} - {\bf \bar{u}_t}$ respectively, evolve as follows:
\begin{align}
    \delta {\bf x_{t+1}} &= A_t \delta {\bf x_t} + B_t \delta {\bf u_t} + S_t(\delta {\bf x_t}) + R_t({\delta {\bf u_t}}) + \epsilon \omega_t, \label{eq:state_perturbation}
    \\
    \delta {\bf u_t} &= K_t \delta {\bf x_t} + \tilde{S}_t(\delta {\bf x_t}) \label{eq:control_perturbation},
\end{align}
where $A_t = \frac{\partial f}{\partial {\bf x}}|_{{(\bf \bar{x}_t}, {\bf \bar{u}_t})}$, $B_t = \frac{\partial f}{\partial {\bf u}}|_{({\bf \bar{x}_t}, {\bf \bar{u}_t})}$, $K_t = \frac{\partial \pi({\bf x_t})}{\partial x}|_{\bf \bar{x}_t}$ and $S_t(.)$, $\tilde{S}_t(.)$, $R_t(.)$ are higher order terms in respective Taylor's expansions. 
By substituting \eqref{eq:control_perturbation} in \eqref{eq:state_perturbation}, we have 
\begin{align}
    \delta {\bf x_{t+1}} = \underbrace{(A_t + B_t K_t)}_{\bar{A}_t} \delta {\bf x_t} + \bar{S}_t(\delta {\bf x_t}) + \epsilon \omega_t,\label{eq:full_state_perturbation}
\end{align}
where $\bar{S}_t(\delta {\bf x_t}) = B_t \tilde{S}_t (\delta {\bf x_t}) + S_t(\delta {\bf x_t}) + R_t(K_t\delta {\bf x_t} + \tilde{S}_t(\delta {\bf x_t}))$. Now, similarly expanding the cost function $c_t({\bf x_t}, {\bf u_t})$ under the policy $ \pi_t({\bf x_t})$ $\forall t\in\{0,1,\cdots, N-1\}$ about the nominal trajectory $\mathbb{T}_{nom}$ as $c_t({\bf x_t}, \pi_t({\bf x_t})) = \bar{c}_t + \bar{C}_t \delta {\bf x_t} + \bar{H}_t(\delta {\bf x_t})$, where $\bar{c}_t = c_t({\bf \bar{x}_t}, \pi_t({\bf \bar{x}_t}))$ and $\bar{C}_t = \frac{\partial c_t({\bf x_t}, \pi_t({\bf x_t}))}{\partial {\bf x_t}}|_{\bf \bar{x}_t}$ results in
\begin{equation}
    J^{\pi} = \sum_{t=0}^{N} \bar{c}_t + \sum_{t=0}^{N} \bar{C}_t \delta {\bf x_t} + \sum_{t=0}^{N} \bar{H}_t(\delta {\bf x_t}).
    \label{eq:cost_to_go_taylors}
\end{equation}

{\bf Lemma 1:} The state perturbation equation \eqref{eq:full_state_perturbation} can be equivalently separated into $\delta {\bf x_t} = \delta{ \bf x_t^l} + \bar{S}_t$ and $\delta {\bf x_{t+1}^l} = \bar{A}_t \delta {\bf x_t^l} + \epsilon \omega_t$, and $\bar{S}_t$ is $O(\epsilon^2)$.\\
{\bf Proof}: Proof is provided in the supplementary file. 

From \eqref{eq:cost_to_go_taylors} and lemma-1, we obtain
\begin{equation}
    J^{\pi} = \underbrace{\sum_{t=0}^{N} \bar{c}_t}_{\bar{J}^{\pi}} + \underbrace{\sum_{t=0}^{N} \bar{C}_t \delta {\bf x_t^l}}_{\delta J^{\pi}_1} + \underbrace{\sum_{t=0}^{N} \bar{H}_t(\delta {\bf x_t}) + \bar{C}_t \bar{S}_t}_{\delta J^{\pi}_2}.
\end{equation}

{\bf Proposition 1:} Given a feedback policy $\pi$ satisfying the smoothness requirements,
\begin{align}
        \tilde{J}^{\pi} &= \mathbb{E}[J^{\pi}] = \bar{J}^{\pi} + O(\epsilon^2),\\
        Var(J^{\pi}) &= \underbrace{Var(\delta J_1^{\pi})}_{O(\epsilon^2)} + O(\epsilon^4).
\end{align}
{\bf Proof:} Proof is provided in the supplementary file.
Noting that ${\bf u_t} = \pi_t({\bf x_t}) = {\bf \bar{u}_t} + K_t \delta {\bf x_t} + \tilde{S}_t (\delta {\bf x_t})$, the above proposition entails the following observations:\\
{\bf i) Observation 1:} The cost-to-go along the nominal trajectory $\mathbb{T}_{nom} = \{{\bf \bar{x}_{0:N}}, {\bf \bar{u}_{0:N-1}}\}$ given by $\bar{J}^{\pi}$ is within  $O(\epsilon^2)$ of the expected cost-to-go $\tilde{J}^{\pi}$ of the policy $\pi$.\\
{\bf ii) Observation 2:} Given the nominal control sequence ${\bf \bar{u}_{0:N-1}}$, the variance of the cost-to-go is overwhelmingly determined by the linear feedback law $K_t \delta {\bf x_t}$ $i.e,$ by the variance of the linear perturbation of the cost-to-go, $\delta J_1^{\pi}$, under the linear closed loop system dynamics $\delta {\bf x_{t+1}} = \bar{A}_t \delta {\bf x_t^l} + \epsilon \omega_t$.

Proposition 1 and the above observations suggest that an optimal open-loop control sequence with a suitable linear feedback law (for the perturbed linear system) wrapped around it is approximately optimal. The following subsection summarizes the problems to be solved for each of them.
\subsection{Decoupled Policy for Feedback Control}
\label{OL_CL_problems}
{\bf Open-Loop Trajectory Design:} The open-loop trajectory is designed by solving the noiseless equivalent of the stochastic optimal control problem \eqref{soc_problem}, as shown below:
\begin{align*}
    \underset{\mathbb{T}_{nom}}{\textnormal{min}} \sum^{N-1}_{1} c_t({\bf \bar{x}_t}, {\bf \bar{u}_t}) + c_N({\bf \bar{x}_N}),\\
    \textnormal{s.t.}~ {\bf \bar{x}_t} = f({\bf \bar{x}_t}, {\bf \bar{u}_t}).
\end{align*}
{\bf Feedback Law:}
A linear feedback law is determined by solving for the optimal feedback gains $K_{1:N-1}^*$ that minimize the variance corresponding to the linear perturbation system around the optimal nominal trajectory $\mathbb{T}^*_{nominal}$ found from the above, as follows:
\begin{align*}
    \underset{\{K_{1:N-1}\}}{\textnormal{min}}~ &Var(\delta J_1^{\pi}),\\
    \delta J_1^{\pi} &= \sum_{t=1}^{N} \bar{C}_t \delta {\bf x_t^l},\\
    \delta {\bf x_{t+1}} &= (A_t + B_t K_t)\delta {\bf x_t^l} + \epsilon \omega_t.
\end{align*}

Albeit convex, the above problem does not have a standard solution. However, since we are really only interested in a good variance rather than the optimal variance, we can solve an inexpensive  time-varying LQR problem as a surrogate to reduce the variance of the cost-to-go. The corresponding problem can be posed as follows:
\begin{align*}
    \underset{{\delta {\bf u_{1:N}}}}{\textnormal{min}}~ \mathbb{E}[\sum_{t=1}^{N-1} \delta {\bf x_t}^T Q_t {\bf x_t} + \delta {\bf u_t}^T R_t {\bf u_t} + \delta {\bf x_N}^T Q_N {\bf x_N}],\\
        \textnormal{s.t.}~\delta {\bf x_{t+1}} = A_t \delta {\bf x_t} + B_t {\bf u_t} + \epsilon \omega_t. 
\end{align*}
\begin{remark}
In fact, if the feedback gain is designed carefully rather than the heuristic LQR above, then one can obtain $O(\epsilon^4)$ near-optimality. This requires some more developments, however, we leave this result out of this paper because of the paucity of space.
\end{remark}
\section{ALGORITHM}
\label{algorithm}

This section presents our decoupled data-based control (D2C-2.0) algorithm. There are two main steps involved in solving for the requisite policy:
\begin{enumerate}
    \item Design an open-loop nominal trajectory by solving the optimization problem via a model-free ILQR based approach $i.e,$ without explicitly making use of an analytic model of the system. This is described in subsection \ref{open_loop_traj_design}. This is the primary difference from the prequel D2C where we used a first order gradient based technique to solve the optimal open loop control problem. 
    \item Determine the parameters of the time-varying linear perturbation system about its optimal nominal trajectory and design an LQR controller based on the system identification. This step is detailed in subsection \ref{lqr_design}.
\end{enumerate}
\subsection{Open-Loop Trajectory Design}
\label{open_loop_traj_design}
An open-loop trajectory is designed by solving the corresponding problem in subsection \ref{OL_CL_problems}. However, since we are solving it in a model-free manner, we need a suitable formulation that can efficiently compute the solution.\\
In this subsection, we present an ILQR based model-free method to solve the open-loop optimization problem. The advantage with ILQR is that the equations involved in it are explicit in system dynamics and their gradients. Hence, in order to make it a model-free algorithm, it is sufficient if we could explicitly obtain the estimates of Jacobians. A sample-efficient way of doing it is described in the following subsubsection. Since ILQR/DDP is a well-established framework, we skip the details and instead present the essential equations in algorithm \ref{model_free_DDP_OL}, algorithm \ref{model_free_DDP_OL_FP} and algorithm \ref{model_free_DDP_OL_BP}, where we also reflect the choice of our regularization scheme.
\subsubsection{Estimation of Jacobians in a Model-Free Setting}
\label{sys_id_solve}

From the Taylor's expansions of `\textit{f}' about the nominal trajectory $({\bf \bar{x}_t}, {\bf \bar{u}_t})$ on both the positive and the negative sides, we obtain the following central difference equation:
\begin{equation}
\begin{split}
    f({\bf \bar{x}_t} + \delta {\bf x_t}, {\bf \bar{u}_t} + \delta {\bf u_t}) - f({\bf \bar{x}_t} - \delta {\bf x_t}, {\bf \bar{u}_t} - \delta {\bf u_t}) =~&2 \begin{bmatrix} f_{\bf x_t} & f_{\bf u_t} \end{bmatrix} \begin{bmatrix}  \delta {\bf x_t} \\ \delta {\bf u_t} \end{bmatrix} +\\ &O(\| \delta {\bf x_t}\|^3 + \| \delta {\bf u_t}\|^3)
     \label{eq : DDP_jacobian_CD}
\end{split}
\end{equation}
Multiplying by $\begin{bmatrix} \delta {\bf x_t}^T & \delta {\bf u_t}^T \end{bmatrix}$ on both sides to the above equation:
\begin{equation*}
    \begin{split}
    \begin{bmatrix}f({\bf \bar{x}_t} + \delta {\bf x_t}, {\bf \bar{u}_t} + \delta {\bf u_t}) - f({\bf \bar{x}_t} - \delta {\bf x_t}, {\bf \bar{u}_t} - \delta {\bf u_t})\end{bmatrix} \times \begin{bmatrix} \delta {\bf x_t}^T & \delta {\bf u_t}^T \end{bmatrix} 
     \\
     =2 \begin{bmatrix} f_{\bf x_t} & f_{\bf u_t} \end{bmatrix} \begin{bmatrix}  \delta {\bf x_t} \delta {\bf x_t}^T & \delta {\bf x_t} \delta {\bf u_t}^T \\ \delta {\bf u_t} \delta {\bf x_t}^T  & \delta {\bf u_t} \delta {\bf u_t}^T \end{bmatrix}+ O(\| \delta {\bf x_t}\|^4 + \| \delta {\bf u_t}\|^4)
     \end{split}
\end{equation*}
Assuming that $\begin{bmatrix}  \delta {\bf x_t} \delta {\bf x_t}^T & \delta {\bf x_t} \delta {\bf u_t}^T \\ \delta {\bf u_t} \delta {\bf x_t}^T  & \delta {\bf u_t} \delta {\bf u_t}^T \end{bmatrix}$ is invertible (which will later be proved to be true by some assumptions), let us perform inversions on either sides of the above equation as follows:
\begin{equation}
    \begin{split}
    \begin{bmatrix} f_{\bf x_t} & f_{\bf u_t} \end{bmatrix} =~&\frac{1}{2}\begin{bmatrix} f({\bf \bar{x}_t} + \delta {\bf x_t}, {\bf \bar{u}_t} + \delta {\bf u_t}) - f({\bf \bar{x}_t} - \delta {\bf x_t}, {\bf \bar{u}_t} - \delta {\bf u_t})\end{bmatrix} \times\\&\begin{bmatrix} \delta {\bf x_t}^T & \delta {\bf u_t}^T \end{bmatrix} \begin{bmatrix}  \delta {\bf x_t} \delta {\bf x_t}^T & \delta {\bf x_t} \delta {\bf u_t}^T \\ \delta {\bf u_t} \delta {\bf x_t}^T  & \delta {\bf u_t} \delta {\bf u_t}^T \end{bmatrix}^{-1} + O(\| \delta {\bf x_t}\|^2 + \| \delta {\bf u_t}\|^2)
     \label{eq : model_free_DDP_jacobian_raw}
     \end{split}
\end{equation}
Equation \eqref{eq : model_free_DDP_jacobian_raw} solves for Jacobians, $f_{\bf x_t}$ and $f_{\bf u_t}$, simultaneously. It is noted that the above formulation requires only 2 evaluations of $f(.)$, given the nominal state and control $\-$ (${\bf \bar{x}_t},{\bf \bar{u}_t}$). The remaining terms (also, the error in the evaluation) are of the order that is quadratic in $\delta {\bf x_t}$ and $\delta {\bf u_t}$. The following extends equation \eqref{eq : model_free_DDP_jacobian_raw} to be used in practical implementations (sampling).\\
We are free to choose the distribution of $\delta {\bf x_t}$ and $\delta {\bf u_t}$. We assume both are i.i.d. Gaussian distributed random variables with zero mean and a standard deviation of $\sigma$.~This ensures that $\begin{bmatrix}  \delta {\bf x_t} \delta {\bf x_t}^T & \delta {\bf x_t} \delta {\bf u_t}^T \\ \delta {\bf u_t} \delta {\bf x_t}^T  & \delta {\bf u_t} \delta {\bf u_t}^T \end{bmatrix}$ is invertible. More on the advantage of using this distribution will be elaborated in the next paragraph. Let $`n_s'$ be the number of samples for each of the random variables, $\delta {\bf x_t}$ and $\delta {\bf u_t}$, as $\delta {\bf X_t} = \begin{bmatrix} \delta {\bf x_t^1}& \delta {\bf x_t^2}& \ldots &\delta {\bf x_t^{n_s}}\end{bmatrix}$ and $\delta {\bf U_t} = \begin{bmatrix} \delta {\bf u_t^1} &\delta {\bf u_t^2}& \ldots& \delta {\bf u_t^{n_s}}\end{bmatrix}$, respectively. Then $\begin{bmatrix} f_{\bf x_t} & f_{\bf u_t} \end{bmatrix}$ is given by the following :
\begin{equation}
    \begin{split}
    \begin{bmatrix} f_{\bf x_t} & f_{\bf u_t} \end{bmatrix} =&\begin{bmatrix} f({\bf \bar{x}_t} + \delta {\bf x_t^1}, {\bf \bar{u}_t} + \delta {\bf u_t^1}) - f({\bf \bar{x}_t} - \delta {\bf x_t^1}, {\bf \bar{u}_t} - \delta {\bf u_t^1})\\ f({\bf \bar{x}_t} + \delta {\bf x_t^2}, {\bf \bar{u}_t} + \delta {\bf u_t^2}) - f({\bf \bar{x}_t} - \delta {\bf x_t^2}, {\bf \bar{u}_t} - \delta {\bf u_t^2}) \\ \vdots  \\f({\bf \bar{x}_t} + \delta {\bf x_t^{n_s}}, {\bf \bar{u}_t} + \delta {\bf u_t^{n_s}}) - f({\bf \bar{x}_t} - \delta {\bf x_t^{n_s}}, {\bf \bar{u}_t} - \delta {\bf u_t^{n_s}})\end{bmatrix} \times\\& \begin{bmatrix} \delta {\bf X_t}^T & \delta {\bf U_t}^T \end{bmatrix} \begin{bmatrix}  \delta {\bf X_t} \delta {\bf X_t}^T & \delta {\bf X_t} \delta {\bf U_t}^T \\ \delta {\bf U_t} \delta {\bf X_t}^T  & \delta {\bf U_t} \delta {\bf U_t}^T \end{bmatrix}^{-1} 
     \label{eq : model_free_DDP_jacobian}
     \end{split}
\end{equation}
Let us consider the terms in the matrix $\delta {\bf XU_t}=\begin{bmatrix}  \delta {\bf X_t} \delta {\bf X_t}^T & \delta {\bf X_t} \delta {\bf U_t}^T \\ \delta {\bf U_t} \delta {\bf X_t}^T  & \delta {\bf U_t} \delta {\bf U_t}^T \end{bmatrix}$.~$\delta {\bf X_t} \delta {\bf X_t}^T = \sum_{i=1}^{n_s} \delta {\bf x_t}^i {\delta {\bf x_t}^i}^T$. Similarly, $\delta {\bf U_t} \delta {\bf U_t}^T = \sum_{i=1}^{n_s} \delta {\bf u_t}^i {\delta {\bf u_t}^i}^T$, $\delta {\bf U_t} \delta {\bf X_t}^T = \sum_{i=1}^{n_s} \delta {\bf u_t}^i {\delta {\bf x_t}^i}^T$ and $\delta {\bf X_t} \delta {\bf U_t}^T = \sum_{i=1}^{n_s} \delta {\bf x_t}^i {\delta {\bf u_t}^i}^T$. From the definition of sample variance, we can write the above matrix as 
\begin{equation*}
\begin{split}
    \delta {\bf XU_t} &= \begin{bmatrix} \sum_{i=1}^{n_s} \delta {\bf x_t}^i {\delta {\bf x_t}^i}^T & \sum_{i=1}^{n_s} \delta {\bf x_t}^i {\delta {\bf u_t}^i}^T \\ \sum_{i=1}^{n_s} \delta {\bf u_t}^i {\delta {\bf x_t}^i}^T & \sum_{i=1}^{n_s} \delta {\bf u_t}^i {\delta {\bf u_t}^i}^T
    \end{bmatrix}\\ &\approx \begin{bmatrix} \sigma^2(n_s - 1) {\text I_{n_x}} & {\text 0_{n_x \times n_u}} \\ 0_{n_u \times n_x} & \sigma^2 (n_s - 1) {\text I_{n_u}}\end{bmatrix} \\
    &= \sigma^2 (n_s - 1){\text I}_{(n_x+n_u) \times (n_x+n_u)}
\end{split}
\end{equation*}
Given that we have high enough number of samples `$n_s$' (typically, slightly more than `$n_x + n_u$' are sufficient), the above approximation holds good. Since the inversion of an identity matrix is trivial and always exists, the above matrix is invertible in equation \eqref{eq : model_free_DDP_jacobian}. Thus, one can calculate $f_{\bf x}$ and $f_{\bf u}$ this way during the backward pass. For the sake of convenience in the later sections, let us refer to this method as `Linear Least Squares by Central Difference (LLS-CD)'. The entire algorithm to determine the optimal nominal trajectory in a model-free fashion is summarized together in Algorithm \ref{model_free_DDP_OL}, Algorithm \ref{model_free_DDP_OL_BP} and Algorithm \ref{model_free_DDP_OL_FP}.
\begin{algorithm}
  \caption{\strut Open-loop trajectory optimization via model-free DDP}
  {\bf Input:} Initial State - ${\bf x_0}$, System parameters - $\mathcal{P}$\;
  $k \gets 1$. ~~\CommentSty{/* Initialize the iteration number $k$ to 1.*/}\\
  $forward\_pass\_flag$ = true. \\
  \CommentSty{/* Run until the difference in costs between subsequent iterations is less an $\epsilon$ fraction of the former cost.*/}\\
  \While {$k==1$ {\textnormal{or}} $({\text cost}(\mathbb{T}_{nom}^{k})/{\text cost}(\mathbb{T}_{nom}^{k-1})) < 1 + \epsilon$}{
  \CommentSty{/*Each iteration has a backward pass followed by a forward pass.*/}\\
  \{$k^{k}_{0:N-1}, K^{k}_{0:N-1}$\}, $backward\_pass\_success\_flag$ $=$ Backward Pass($\mathbb{T}_{nom}^{k}$, $\mathcal{P}$).\\
  \If{backward\_pass\_success\_flag == true}{
      $\mathbb{T}_{nom}^{k+1}, forward\_pass\_flag$ $=$ Forward Pass($\mathbb{T}_{nom}^{k}$,$\{k^{k}_{0:N-1}, K^k_{0:N-1}\}$, $\mathcal{P}$).\\
      \While{forward\_pass\_flag == false}{
            $\mathbb{T}_{nom}^{k+1}, forward\_pass\_flag$ $=$ Forward Pass($\mathbb{T}_{nom}^{k}$,$\{k^{k}_{0:N-1}, K^k_{0:N-1}\}$, $\mathcal{P}$).\\
            Reduce $\alpha$ from $\mathcal{P}$.    
        }
      }
   \Else{
        Increase $\mu$ from $\mathcal{P}$.~~\CommentSty{/* Regularization step */}\\
    }
      $k \leftarrow k + 1$. \\
  $\mathbb{T}^{*}_{nom} \gets \mathbb{T}^{k+1}_{nom}$.\\
  }
  
  {\bf return $\mathbb{T}^{*}_{nom}$}\\
  \label{model_free_DDP_OL}
\end{algorithm}
\begin{algorithm}
  \caption{\strut Forward Pass}
  {\bf Input:} Nominal trajectory - $\mathbb{T}_{nom}^{k}$, previous iteration policy parameters - $\{k_{0:N-1}, K_{0:N-1}\}$ and system and cost parameters - $\mathcal{P}$.\
  $\{{\bf \bar{x}_t^{prev}, \bar{u}_t^{prev}}\} \gets \mathbb{T}_{nom}^{k}$.\\
  $t \gets 0$.\\
  \While {$t < N$}{
    \CommentSty{/*Simulate one step forward ($\alpha$ is the line-search parameter.)*/}
    \begin{equation*}
    \begin{split}
        {\bf \bar{u}_{t}} &= {\bf \bar{u}_t^{prev}} + \alpha k_t + K_{t} ({\bf \bar{x}_t} - {\bf \bar{x}_t^{prev}}),\\
        {\bf \bar{x}_{t+1}} &= simulate\_forward\_step({\bf \bar{x}_t}, {\bf \bar{u}_t}).
        \end{split}
    \end{equation*}
    $t \leftarrow t + 1$. 
    }
  $\mathbb{T}_{nom}^{k+1}$ $\gets$ $\{{\bf \bar{x}_{0:N}}, {\bf \bar{u}_{0:N-1}}\}.$\\
  \If{$\mathbb{T}_{nom}^{k+1}$ to $\mathbb{T}_{nom}^{k}$ \textnormal{is acceptable}}{
    {\bf return} $\mathbb{T}_{nom}^{k+1}$, true.
  }
  \Else{
    {\bf return} $\mathbb{T}_{nom}^{k}$, false.
  }
  \label{model_free_DDP_OL_FP}
\end{algorithm}
\begin{algorithm}
  \caption{\strut Backward Pass}
  {\bf Input:} Nominal trajectory - $\mathbb{T}_{nom}^{k}$, previous iteration policy parameters - $\{k_{0:N-1}, K_{0:N-1}\}$, horizon - N and system and cost parameters - $\mathcal{P}$.\\
  \CommentSty{/* Backward pass starts from the final time-step i.e, N-1.*/}\\
  $t \gets N - 1$.\ \\
  Compute $J_{\bf x_N}$ and $J_{\bf x_N x_N}$ using boundary conditions.\\
  \CommentSty{/*Keep a copy of previous policy gains.*/}\\
  $k\_old \gets k_{0:N}$ {and}  $K\_old \gets K_{0:N}$.\\
  
  \While {$t >= 0$}{
  \CommentSty{/*Obtain the Jacobians from simulator rollouts as shown in equation \eqref{eq : model_free_DDP_jacobian}:*/}\\
  $f_{\bf x_t}, f_{\bf u_t} \gets model\_free\_jacobian({\bf \bar{x}_t},{\bf \bar{u}_t}).$\\
  \CommentSty{/*Obtain the partials of the Q function as follows:*/}
    \begin{equation*}
    \begin{split}
    Q_{\bf x_t} &= c_{\bf x_t}  + f_{\bf x_t}^T J_{\bf x_{t+1}}^{\prime},\\
    Q_{\bf u_t} &= c_{\bf u_t} + f_{\bf u_t}^T J_{\bf x_{t+1}}^{\prime},\\
    Q_{\bf x_t x_t} &= c_{\bf x_t x_t} + f_{\bf x_t}^T J_{\bf x_{t+1} x_{t+1}}^{\prime} f_{\bf x_t}, \\
    Q_{\bf u_t x_t} &= c_{\bf u_t x_t} + f_{\bf u_t}^T (J_{\bf x_{t+1} x_{t+1}}^{\prime} + \mu I_{n_x \times n_x}) f_{\bf x_t}, \\
    Q_{\bf u_t u_t} &= c_{\bf u_t u_t} + f_{\bf u_t}^T (J_{\bf x_{t+1} x_{t+1}}^{\prime} + \mu I_{n_x \times n_x}) f_{\bf u_t}.\\
    \end{split}
\end{equation*}
  \If{$Q_{\bf u_t u_t}$ {\textnormal{is positive-definite}}}{
    \begin{equation*}
    \begin{split}
        k_t &= -Q_{\bf u_t u_t}^{-1} Q_{\bf u_t},\\
        K_t &= -Q_{\bf u_t u_t}^{-1} Q_{\bf u_t x_t}.
        \end{split}
    \end{equation*}
    }
    
  \Else{
    \CommentSty{/*If $Q_{\bf u_t u_t}$ is not positive-definite, then, abort the backward pass.*/}\\
    {\bf return } $\{ k\_old, K\_old\}$, false.
  }
  \CommentSty{/*Obtain the partials of the value function $J_t$ as follows:*/}
  \begin{equation*}
      \begin{split}
          J_{\bf x_t} &= Q_{\bf x_t} + K_{t}^T Q_{\bf u_t u_t} k_t + K_t^T Q_{\bf u_t} + Q_{\bf u_t x_t}^T k_t,\\
          J_{\bf x_t x_t} &= Q_{\bf x_t x_t} + K_t^T Q_{\bf u_t u_t} K_t + K_t^T Q_{\bf u_t x_t} + Q_{\bf u_t x_t}^T K_t.
      \end{split}
  \end{equation*}
  $t \leftarrow t - 1$ 
  }
  $k\_new = k_{0:N-1},$\\
  $K\_new = K_{0:N-1}.$\\
  {\bf return} $\{k\_new, K\_new\}$, true.
  \label{model_free_DDP_OL_BP}
\end{algorithm}
\subsection{Feedback Law}
\label{lqr_design}
As discussed in the previous section, a time-varying LQR feedback law is used as a surrogate for the minimum variance problem. The LQR gains $\{K_{1:N-1}\}$ can be fully determined if $Q_t$, $R_t$, $A_t$ and $B_t$ are known. The first two are the tuning parameters in the cost function. Determining $A_t$ and $B_t$ for a linear system is a system identification problem that can be solved by `LLS-CD' as described under subsubsection \ref{sys_id_solve}. 
\section{SIMULATIONS}
\label{simulations}
This section reports the implementation details of D2C-2.0 in simulation on various systems and their training results. We train the systems with varying complexities in terms of the dimension of the state and the action spaces. The physical models of the system are deployed in the simulation platform `MuJoCo-2.0' \cite{mujoco} as a surrogate to their analytical models. The models are imported from the OpenAI gym \cite{openai gym} and Deepmind's control suite \cite{control suite}.

\begin{figure}
\centering
\begin{multicols}{2}
    \begin{subfigure}{0.8\linewidth}
    {\includegraphics[width=\linewidth, height=3.8cm]{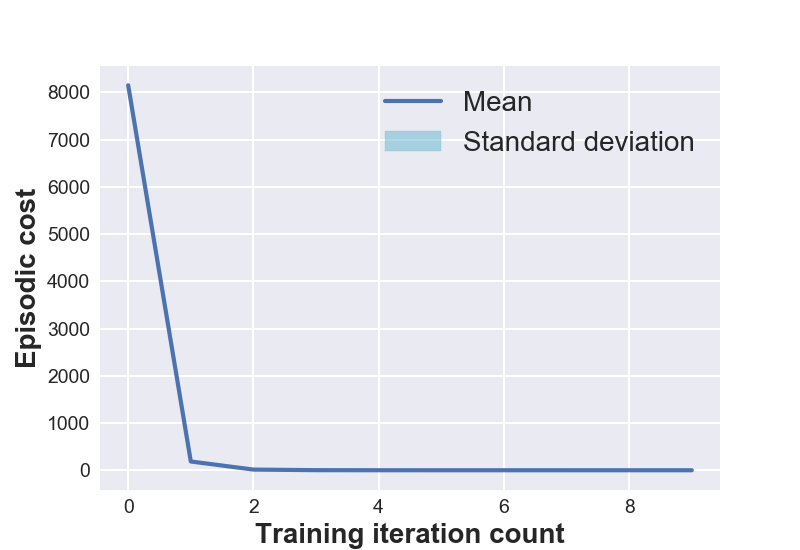}}
    \caption{Inverted Pendulum}
    \end{subfigure}
    \begin{subfigure}{0.82\linewidth}
        {\includegraphics[width=\linewidth, height=3.8cm]{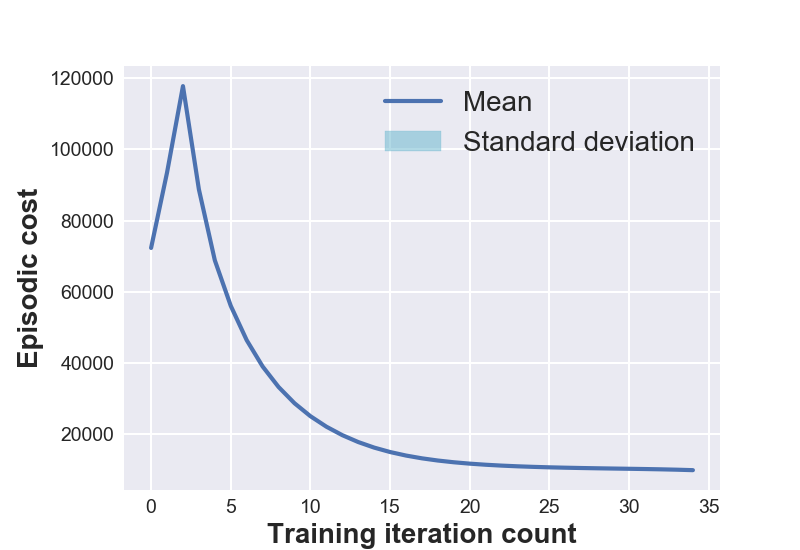}}
    \caption{Cart-Pole}
    \end{subfigure}
     \begin{subfigure}{0.8\linewidth}
        {\includegraphics[width=\linewidth, height=3.8cm]{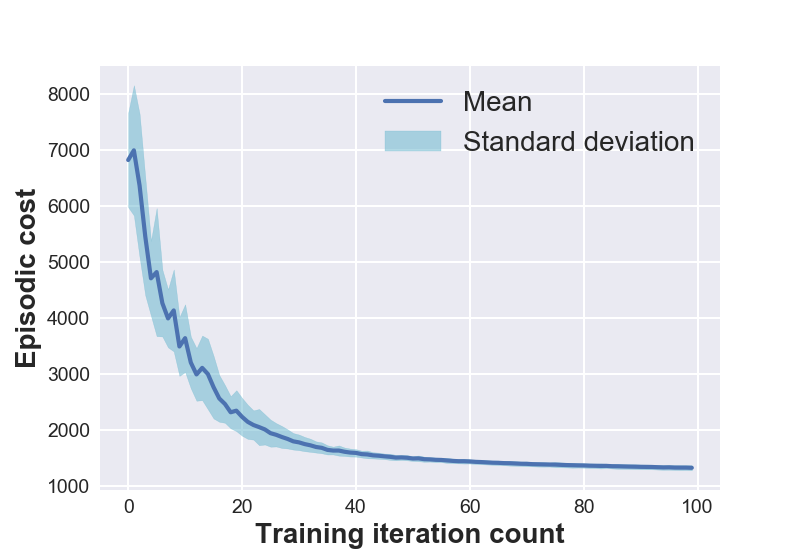}}
    \caption{3-link Swimmer}
    \end{subfigure}
      \begin{subfigure}{0.8\linewidth}
        {\includegraphics[width=\linewidth, height=3.8cm]{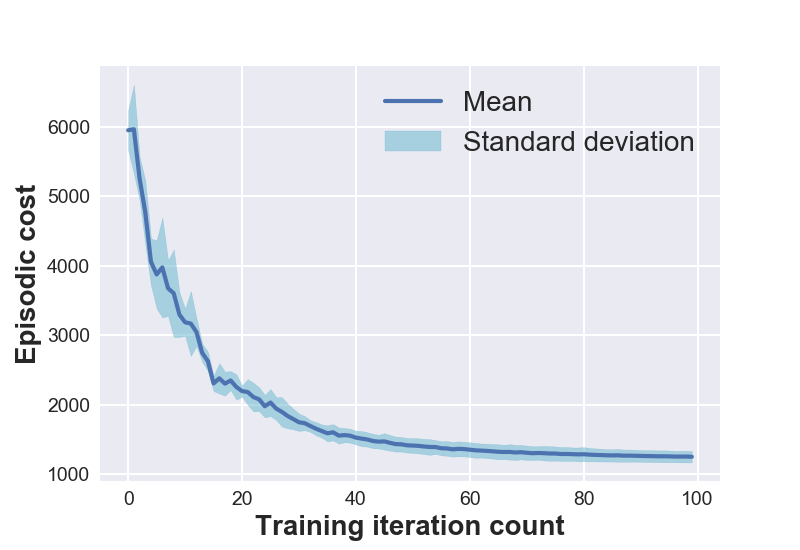}}
      \caption{6-link swimmer}
      \end{subfigure}
      \begin{subfigure}{0.8\linewidth}
        {\includegraphics[width=\linewidth, height=3.8cm]{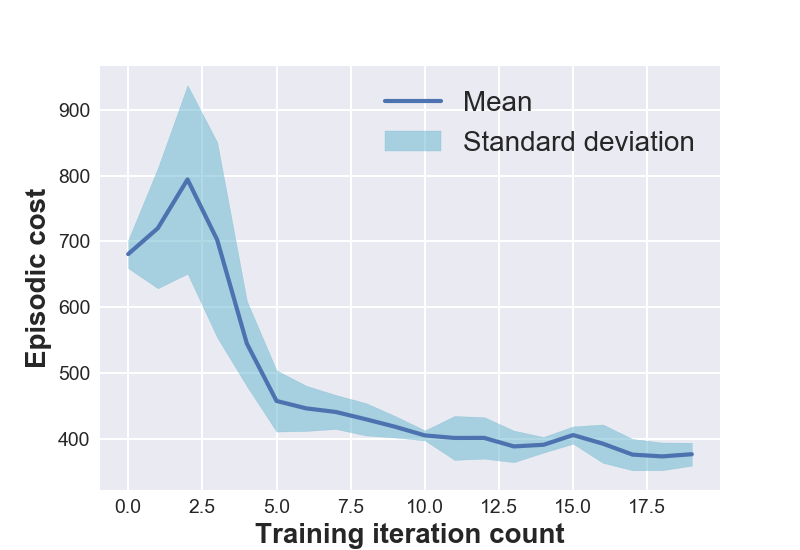}}
        \caption{Robotic fish}
      \end{subfigure}
        \begin{subfigure}{0.8\linewidth}
    {\includegraphics[width=\linewidth, height=3.8cm]{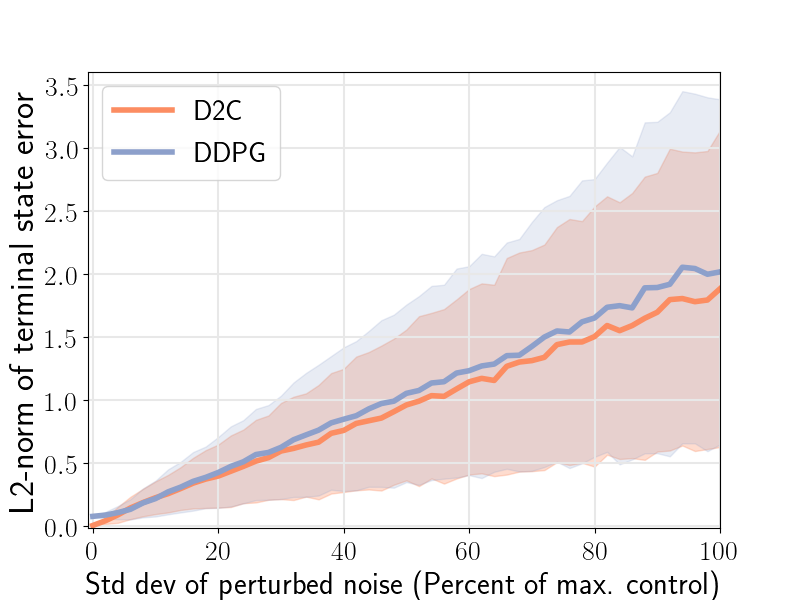}}
    \caption{Inverted Pendulum}
    \end{subfigure}
    \begin{subfigure}{0.8\linewidth}
        {\includegraphics[width=\linewidth, height=3.8cm]{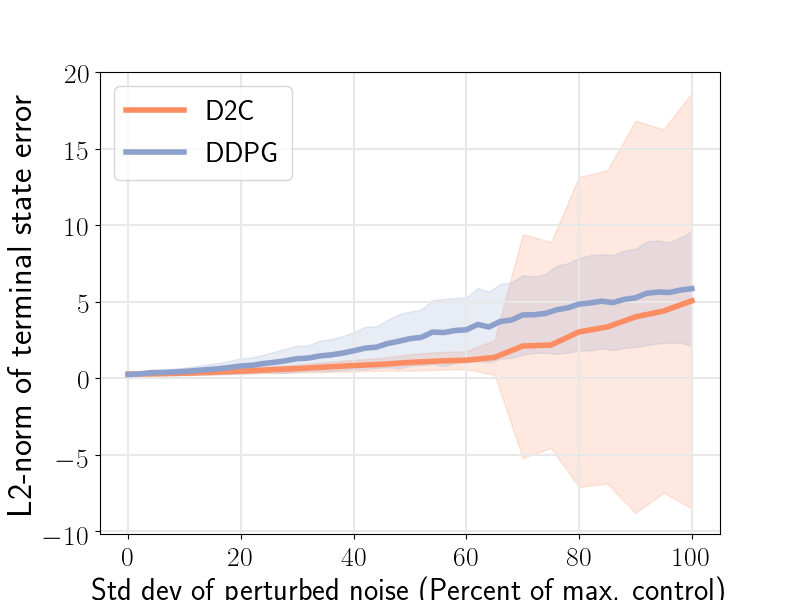}}
    \caption{Cart-Pole}
    \end{subfigure}
          \begin{subfigure}{0.8\linewidth}
        {\includegraphics[width=\linewidth, height=3.8cm]{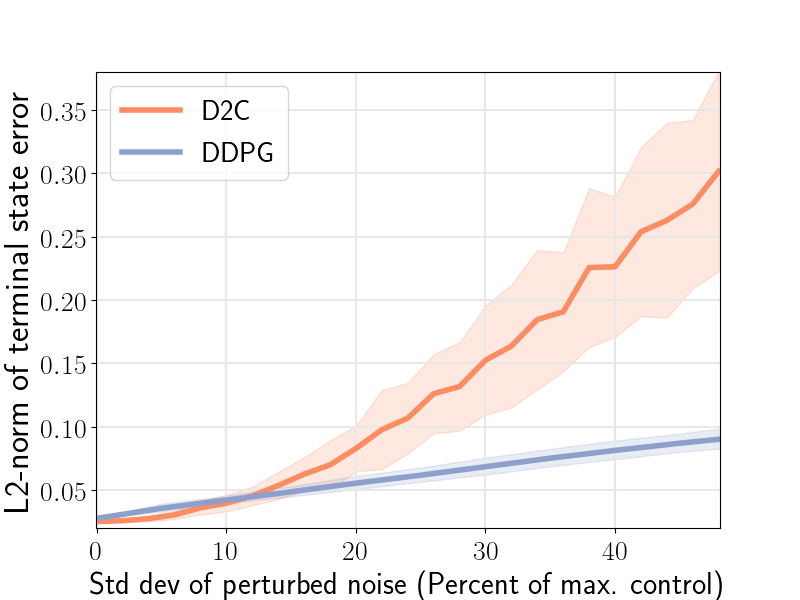}}
    \caption{3-link swimmer}
    \end{subfigure}
      \begin{subfigure}{0.8\linewidth}
        {\includegraphics[width=\linewidth, height=3.8cm]{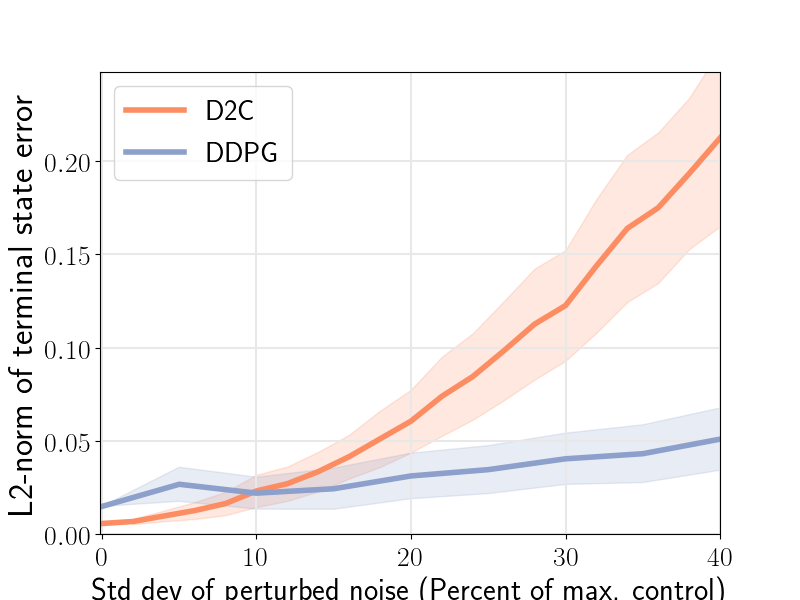}}
      \caption{6-link swimmer}
      \end{subfigure}
      \begin{subfigure}{0.8\linewidth}
        {\includegraphics[width=\linewidth, height=3.8cm]{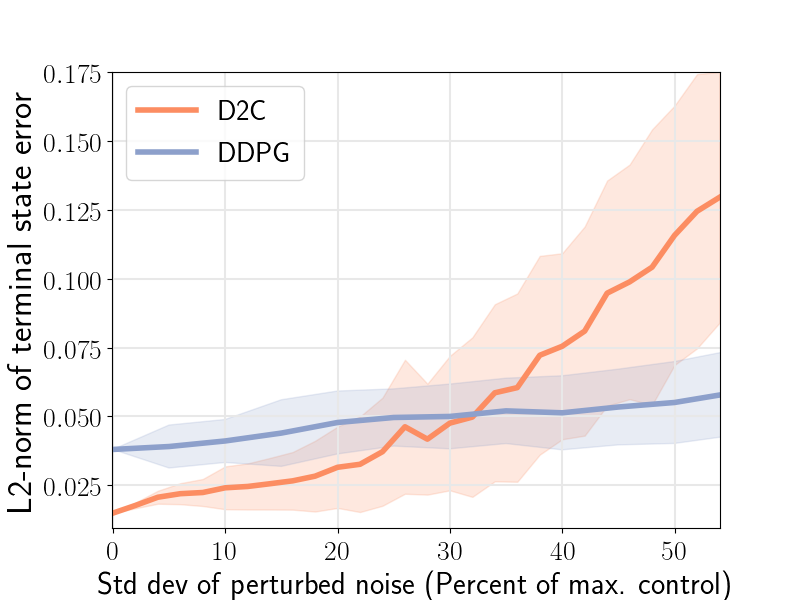}}
      \caption{Fish}
      \end{subfigure}
\end{multicols}
\caption{(Left) Episodic cost vs.~training iteration number in D2C-2.0 (Right) Terminal state MSE during testing in D2C (2.0) vs DDPG}
\label{d2c_2_training_testing}
\end{figure}

Figure \ref{d2c_2_training_testing} (left) shows the open-loop training plots of various systems under consideration. The data is averaged over 5 experimental runs to reflect the training variance and the repeatability of D2C-2.0 over multiple training sessions. One might notice that the episodic cost is not monotonically decreasing during training. This is because of the flexibility provided in training for the episodic cost to stay within a `band' around the cost at the previous iteration. Though this might momentarily increase the cost in subsequent training iterations, it is observed to have improved the overall speed of convergence.

\section{DISCUSSION}
\label{discussion}
In this section, we compare D2C-2.0 with D2C and DDPG, and discuss the merits and demerits of each of the methods.

The manuscript on D2C provided an in-depth analysis on the comparison of D2C with DDPG. We have shown that training in D2C is much more sample-efficient and fares well w.r.t. the training time and the reproducibility of the resulting policy.


\begin{table}[h]
\caption{Comparison of the simulation parameters and training outcomes of D2C-2.0 with other baselines.}
\begin{center}
\begin{tabular}{|c|c|c|c|c|}
\hline
& \multicolumn{4}{|c|}{\bf Training time (in sec.)}\\\cline{2-5}
{\bf System}&\multicolumn{2}{|c|}{\bf D2C}& &\\\cline{2-3}
& {\bf Open-}& {\bf Closed-}&{\bf DDPG}&{\bf D2C-}\\
&{\bf loop}&{\bf loop}&&{\bf 2.0}\\
\hline
{\bf Inverted} &12.9  &$<0.1$& 2261.15&0.33\\
{\bf Pendulum} & & & &\\
\hline
{\bf Cart pole}& 15.0 & 1.33& 6306.7 &1.62\\
\hline
{\bf 3-link} & & &  &\\
{\bf Swimmer}& 4001.0&4.6 & 38833.64 & 186.2\\
\hline
{\bf 6-link} &  & & &\\
{\bf Swimmer} &3585.4&16.4 & 88160 & 127.2\\ 
\hline
{\bf Fish} & 6011.2 &75.6 & 124367.6 & 54.8\\
\hline
\end{tabular}
\end{center}
\label{d2c_comparison_table2}
\end{table}
Table \ref{d2c_comparison_table2} shows the training times of D2C, DDPG and D2C-2.0. It is evident from the table that for simple and lower dimensional systems such as pendulum and cart pole, the policy can be almost calculated online (parallelization could make it much faster). It is due to the tendency of ILQR/DDP to quickly converge if the linearization of system models has bigger basin of attraction (such as pendulum, which is, in fact, often approximated to be linear about its resting state). Moreover, for a higher dimensional system such as in the case of a robotic fish, D2C-2.0 took 54.8 seconds while the original D2C took around 6096 seconds, which is approximately 90 times faster. The table makes it evident that D2C-2.0 clearly exhibits much faster convergence when compared to both D2C and DDPG, while D2C is much faster than DDPG. The primary reason for the much larger training times in DDPG is due to the complex parametrization of the feedback policy in terms of deep neural networks which runs into hundreds of thousands of parameters while D2C has a parametrization, primarily the open loop control sequence which is of dimension $pT$, where $p$ is the number of inputs and $T$ is the time horizon of the problem,  that runs from tens to a few thousands (see Table \ref{parasize}).  \\

\begin{table}[ht]
\caption{Parameter size comparison between D2C and DDPG}
\label{parasize}
\centering
\vspace{0.1in}
\begin{threeparttable}
\setlength{\tabcolsep}{0.6mm}{
\begin{tabular}{|c|c|c|c|c|}
\hline
System& No. of  & No. of & No. of  &No. of \\
&steps&actuators&parameters&parameters\\
&&&optimized&optimized\\
&&& in D2C& in DDPG\\
\hline
Inverted&30  &1&30&244002\\
Pendulum& & &&\\
\hline
Cart pole& 30 &1& 30&245602\\
\hline
3-link & 1600&2 & 3200&251103\\
Swimmer& & &&\\
\hline
6-link &1500&5 &7500&258006\\
Swimmer &&&& \\ 
\hline
Fish&1200 &6 & 7200&266806\\
\hline
\end{tabular}
}
\vspace{0.1in}
\end{threeparttable}
\end{table}

A major reason behind the differences in training efficiencies stems from their formulation. D2C is based on gradient descent, which is linearly convergent. It is a direct and generic method of solving an optimization problem. On the other hand, the equations involved in D2C-2.0 arise from the optimality conditions that inherently exploit the recursive optimal structure and the causality in the problem. Moreover, the equations are explicit in the Jacobians of the system model. As a result, it is possible to solve for it {\it independently} and in a sample \-efficient manner, at every time-step of each iteration. 
However, one should expect that as we consider very high dimensional systems, the algorithmic complexity of D2C-2.0 could outgrow that of D2C. This is because it is $O(n_x^2)$  in terms of its computational complexity. Also, D2C, by the nature of its formulation does not need the state of the system, and hence, is generalizable to problems with partial state observability. 

The closed-loop performance of D2C-2.0 w.r.t. that of DDPG is shown in the figure \ref{d2c_2_training_testing} (right). Noise is injected into the system via the control channel by varying the noise scaling factor-$\epsilon$ and then, measuring the performance. We note here that DDPG fails to converge for the robotic fish problem, and thus, the nominal performance ($\epsilon =0$) is poor. We see that the performance of D2C is better than that of DDPG when the noise levels are lower, however, the performance of DDPG remains relatively similar for higher levels of noise while that of D2C degrades. In our opinion, this closed loop performance aspect of the methods needs to be studied in more detail in order to make a definitive judgement, and is left for future work. For instance, by the time DDPG overtakes D2C in terms of the average performance, the actual performance of both systems is quite poor, and thus, it is not clear if this can be construed as an advantage.

Now, as mentioned in section \ref{related_work}, earlier works have attempted to tackle the idea of model-free ILQR and have mostly confined to finite-differences for Jacobian computation. Table \ref{jacobian_comparison_table} shows the comparison of per-iteration computational times between `finite-differences' and `linear least squares with central difference formulation'. It is clearly evident that, as the dimension of the state space increases, the method of finite-differences requires many more function evaluations, and hence, our LLS-CD method is much more efficient.
\begin{table}[h]
\caption{Comparison of the computational times (in sec.) per iteration  (averaged over 5 runs).}
\begin{center}
\begin{tabular}{|c|c|c|}
\hline
{\bf System}&{\bf FD}&{\bf LLS-CD}\\
\hline
{\bf Inverted} &  &\\
{\bf Pendulum} &0.017 &0.033\\
\hline
{\bf Cart pole}& 0.0315 &0.0463 \\
\hline
{\bf Acrobot}& 0.554 &0.625 \\
\hline
{\bf 3-link} & &\\
{\bf Swimmer}& 4.39&1.86 \\
\hline
{\bf 6-link} &  & \\
{\bf Swimmer} &14.43& 1.278  \\ 
\hline
{\bf Fish} & 34.75 & 2.74\\
\hline
\end{tabular}
\end{center}
\label{jacobian_comparison_table}
\end{table}

\section{Conclusion}
In order to conclude the discussion that started the global vs local policy approximation dilemma in section \ref{introduction}, how well deep RL methods extrapolate their policy to the entire continuous state and action spaces remains an open question. Now, consider that the results shown in the Table \ref{d2c_comparison_table2} are based on serial implementations on an off-the-shelf computer and D2C-2.0 can be highly parallelizable. By looking at the order of magnitude of numbers in the same table, we expect with reasonable augmentation of computational power by parallelization, D2C-2.0 could offer a near-real time solution in high dimensional problems. In such cases, one could rely on the near-optimal policy described in this paper for low noises, and re-solve for the open-loop trajectory online by the approach presented in this paper along with the attendant feedback, whenever the noise takes the method out of the region of attraction of the linear feedback policy. This will be our motivation for the future work along with addressing the following limitation. \\
{\bf Limitations:} D2C-2.0, in its current form, is not reliable for hybrid models such as legged robots that involve impacts. Primarily, the issue is in the feedback design which requires  smoothness of the underlying dynamics while legged systems are only piece-wise smooth. Future work shall model impacts to tractably incorporate them in order to estimate the required gradients.  This seems plausible as iLQR has been used to design controllers for such systems in the literature \cite{ILQG_complex_behaviors}.


\section{Supplementary}
{\bf Proof of Lemma 1:}

We proceed by induction. The first general instance of the recursion occurs at $t=3$.
It can be shown that: 
 \begin{align}
 &\delta {\bf x_3} = \underbrace{(\bar{A}_2\bar{A}_1(\epsilon w_0) + \bar{A}_2 (\epsilon w_1) + \epsilon w_2)}_{\delta {\bf x}_3^l} + \nonumber\\
 &\underbrace{\{\bar{A}_2 \bar{S}_1(\epsilon w_0) + \bar{S}_2(\bar{A}_1(\epsilon w_0) + \bar{S}_2(\bar{A}_1(\epsilon w_0) + \epsilon w_1 + \bar{S}_1(\epsilon w_0))\}}_{\bar{\bar{S}}_3}. 
 \end{align}
 Noting that $\bar{S}_1(.)$ and $\bar{S}_2(.)$ are second and higher order terms, it follows that $\bar{\bar{S}}_3$ is $O(\epsilon^2)$.w \\
 Suppose now that $\delta {\bf x_t} = \delta {\bf x_t^l} + \bar{\bar{S}}_t$ where $\bar{\bar{S}}_t$ is $O(\epsilon^2)$. Then:
 \begin{align}
 \delta {\bf x_{t+1}} = \bar{A}_{t+1}(\delta {\bf x_t^l} + \bar{\bar{S_t}}) + \epsilon w_t + \bar{S}_{t+1}(\delta {\bf x_t}), \nonumber\\
 = \underbrace{(\bar{A}_{t+1} \delta {\bf x_t^l} + \epsilon w_t)}_{\delta {\bf x_{t+1}^l}} +\underbrace{\{\bar{A}_{t+1}\bar{\bar{S}}_t + \bar{S}_{t+1}(\delta {\bf x_t})\}}_{\bar{\bar{S}}_{t+1}}.
 \end{align}
 Noting that $\bar{S}_{t+1}$ is $O(\epsilon^2)$ and that $\bar{\bar{S}}_{t+1}$ is $O(\epsilon^2)$ by assumption, the result follows. 

{\bf {Proof of Proposition 1}}:\\

From (6), we get,
\begin{align} 
\tilde{J}^{\pi} = \mathbb{E}[J^{\pi}] =  \mathbb{E}[ \bar{J}^{\pi} + \delta J_1^{\pi} + \delta J_2^{\pi}], \nonumber\\
= \bar{J}^{\pi} + \mathbb{E}[\delta J_2^{\pi}]  = \bar{J}^{\pi} + O(\epsilon^2), \label{eq.10}
\end{align}
The first equality in the last line of the equations before follows from the fact that $\mathbb{E}[\delta {\bf x_t^l}] = 0$, since its the linear part of the state perturbation driven by white noise and by definition $\delta {\bf x_1^l} = 0$. The second equality follows from the fact that $\delta J_2^{\pi}$ is an $O(\epsilon^2)$ function.  Now,
\begin{align}
\text{Var}(J^{\pi}) = \mathbb{E}[ J^{\pi} - \tilde{J}^{\pi}]^2 \nonumber\\
= \mathbb{E}[ \bar{J}_0^{\pi} + \delta J_1^{\pi} + \delta J_2^{\pi} - \bar{J}_0^{\pi} - \delta \tilde{J}_2^{\pi}]^2 \nonumber\\
= \text{Var}(\delta J_1^{\pi}) + \text{Var}(\delta J_2^{\pi})  + 2 \mathbb{E}[\delta J_1^{\pi} \delta J_2^{\pi}].
\end{align}
Since $\delta J_2^{\pi}$ is $O(\epsilon^2)$, $\text{Var}(\delta J_2^{\pi})$ is an $O(\epsilon^4)$ function. It can be shown that $\mathbb{E}[\delta J_1^{\pi} \delta J_2^{\pi}]$ is $O(\epsilon^4)$ as well (proof is given in Lemma 2). Finally $\text{Var}(\delta J_1^{\pi})$ is an $O(\epsilon^2)$ function because $\delta {\bf x^l}$ is an $O(\epsilon)$ function. Combining these, we will get the desired result.

{\bf{Lemma 2}:} \textit{Let $\delta J_1^{\pi}$, $\delta J_2^{\pi}$ be as defined in (6). Then, $\mathbb{E} [\delta J_1 \delta J_2]$ is an $O(\epsilon^4)$ function.}\\
\vspace{3cm}

{\bf Proof of Lemma 2:}\\
In the following, we suppress the explicit dependence on $\pi$ for $\delta J_1^{\pi}$ and $\delta J_2^{\pi}$ for notational convenience.
Recall that $\delta J_1 = \sum_{t=0}^T c_t^x \delta {\bf x_t^l}$, and $\delta J_2 = \sum_{t=0}^T \bar{H}_t(\delta {\bf x_t}) + c_t^x \bar{\bar{S}}_t$.  For notational convenience, let us consider the scalar case, the vector case follows readily at the expense of more elaborate notation. Let us first consider $\bar{\bar{S}}_2$. We have that $\bar{\bar{S}}_2 = \bar{A}_2\bar{S}_1(\epsilon w_0) + \bar{S}_2(\bar{A}_1(\epsilon w_0) + \epsilon w_1+ \bar{S}_1(\epsilon w_0))$. Then, it follows that:
\begin{align}
\bar{\bar{S}}_2 = \bar{A}_2 \bar{S}_1^{(2)}(\epsilon w_0)^2 + \bar{S}_2^{(2)}(\bar{A}_1 \epsilon w_0 + \epsilon w_1)^2 + O(\epsilon^3),
\end{align}
where $\bar{S}_t^{(2)}$ represents the coefficient of the second order term in the expansion of $\bar{S}_t$. A similar observation holds for $H_2(\delta {\bf x_2})$ in that:
\begin{align}
\bar{H}_2(\delta {\bf x_2}) = \bar{H}_2^{(2)}(\bar{A}_1(\epsilon w_0) + \epsilon w_1)^2 + O(\epsilon^3),
\end{align}
where $\bar{H}_t^{(2)}$ is the coefficient of the second order term in the expansion of $\bar{H}_t$. Note that $\epsilon w_0 = \delta {\bf x_1^l}$ and $\bar{A}_1(\epsilon w_0) + \epsilon w_1= \delta {\bf x_2^l}$. Therefore, it follows that we may write:
\begin{align}
\bar{H}_t(\delta {\bf x_t}) + C_t^x \bar{\bar{S}}_t = \sum_{\tau = 0}^{t-1} q_{t,\tau}(\delta {\bf x_{\tau}^l})^2 + O(\epsilon^3),
\end{align}
for suitably defined coefficients $q_{t,\tau}$. 
Therefore, it follows that 
\begin{align}
\delta J_2 = \sum_{t=1}^T \bar{H}_t(\delta {\bf x_t}) + C_t^x \bar{\bar{S}}_t\nonumber\\
= \sum_{\tau = 0}^T \bar{q}_{T,\tau}(\delta {\bf x_{\tau}^l})^2+ O(\epsilon^3),
\end{align}
for suitably defined $\bar{q}_{T,\tau}$. Therefore:
\begin{align}
\delta J_1 \delta J_2 = \sum_{t,\tau} C_{\tau}^x(\delta {\bf x_{\tau}^l})\bar{q}_{T,t}(\delta {\bf x_t^l})^2 + O(\epsilon^4).
\end{align}
Taking expectations on both sides:
\begin{align}
E[\delta J_1 \delta J_2] = \sum_{t,\tau} C_{\tau}^x \bar{q}_{T,t} E[\delta {\bf x_{\tau}^l} (\delta {\bf x_t^l})^2] + O(\epsilon^4).
\end{align}
Break $\delta {\bf x_t^l} = (\delta {\bf x_t^l} - \delta {\bf x_{\tau}^l}) + \delta {\bf x_{\tau}^l}$, assuming $\tau < t$. Then, it follows that:
\begin{align}
E[\delta {\bf x_{\tau}^l} (\delta {\bf x_t^l})^2] = E[\delta {\bf x_{\tau}^l} (\delta {\bf x_t^l} - \delta {\bf x_{\tau}^l})^2] + E[(\delta {\bf x_{\tau}^l)^3}]  \nonumber\\
+ 2 E[(\delta {\bf x_t^l} - \delta {\bf x_{\tau}^l})(\delta {\bf x_{\tau}^l})^2]\nonumber\\
= E[(\delta {\bf x_{\tau}^l})^3],
\end{align}
where the first and last terms in the first equality drop out due to the independence of the increment $(\delta {\bf x_t^l} - \delta {\bf x_{\tau}^l})$ from $\delta {\bf x_{\tau}^l}$, and the fact that $E[\delta {\bf x_t^l} - \delta {\bf x_{\tau}^l}] = 0$ and $E[\delta {\bf x_{\tau}^l}] = 0$. Since $\delta {\bf x_{\tau}^l}$ is the state of the linear system $\delta {\bf x_{t+1}}= \bar{A}_t \delta {\bf x_t^l} + \epsilon w_t$, it may again be shown that:
\begin{align}
E[\delta x_{\tau}^l]^3 = \sum_{s_1,s_2,s_3} \Phi_{\tau, s_1}\Phi_{\tau,s_2}\Phi_{\tau,s_3} E[w_{s_1}w_{s_2}w_{s_3}],
\end{align}
where $\Phi_{t,\tau}$ represents the state transitions operator between times $\tau$ and $t$, and follows from  the closed loop dynamics. Now, due to the independence of the noise terms $w_t$, it follows that $E[w_{s_1}w_{s_2}w_{s_3}] = 0$ regardless of $s_1,s_2,s_3$.\\
 An analogous argument as above can be repeated for the case when $\tau > t$. Therefore, it follows that $E[\delta J_1 \delta J_2] = O(\epsilon^4)$.




%
{\bf DDPG Algorithm:} Deep Deterministic Policy Gradient (DDPG) is a policy-gradient based off-policy reinforcement learning algorithm that operates in continuous state and action spaces. It relies on two function approximation networks one each for the actor and the critic. The critic network estimates the $Q(s,a)$ value given the state and the action taken, while the actor network engenders a policy given the current state. Neural networks are employed to represent the networks. 

The off-policy characteristic of the algorithm employs a separate behavioural policy by introducing additive noise to the policy output obtained from the actor network. The critic network minimizes loss based on the temporal-difference (TD) error and the actor network uses the deterministic policy gradient theorem to update its policy gradient as shown below:

Critic update by minimizing the loss: 
\begin{equation*}
    L = \frac{1}{N} \Sigma_{i}(y_i - Q(s_i, a_i| \theta^{Q}))^2
\end{equation*}

Actor policy gradient update:
\begin{equation*}
    \nabla_{\theta^{\mu}} \approx \frac{1}{N} \Sigma_i \nabla_a Q(s,a|\theta^{Q})|_{s=s_i,a=\mu(s_i)} \nabla_{\theta^{\mu}} \mu(s|\theta^{\mu})|_{s_i}
\end{equation*}

The actor and the critic networks consists of two hidden layers with the first layer containing 400 ('{\it relu}' activated) units followed by the second layer containing 300 ('{\it relu}' activated) units. The output layer of the actor network has the number of ('tanh' activated) units equal to that of the number of actions in the action space.  

Target networks one each for the actor and the critic are employed for a gradual update of network parameters, thereby reducing the oscillations and a better training stability. The target networks are updated at $\tau = 0.001$. Experience replay is another technique that improves the stability of training by training the network with a batch of randomized data samples from its experience. We have used a batch size of 32 for the inverted pendulum and the cartpole examples, whereas it is 64 for the rest. Finally, the networks are compiled using Adams' optimizer with a learning rate of 0.001. 

To account for state-space exploration, the behavioural policy consists of an off-policy term arising from a random process. We obtain discrete samples from Ornstein-Uhlenbeck (OU) process to generate noise as followed in the original DDPG method. The OU process has mean-reverting property and produces temporally correlated noise samples as follows:
\begin{equation*}
    dx_t = \Theta (\mu - x_t)dt + \sigma dW
\end{equation*}
where $\Theta$ indicates how fast the process reverts to mean, $\mu$ is the equilibrium or the mean value and $\sigma$ corresponds to the degree of volatility of the process. $\Theta$ is set to 0.15, $\mu$ to 0 and $\sigma$ is annealed from 0.35 to 0.05 over the training process. 

\end{document}